\documentclass[11pt,a4,fleqn]{article}
\usepackage{graphicx}
\input{amssym}
\newtheorem{theorem}{\bf Theorem}[section]

\newtheorem{lemma}[theorem]{\bf Lemma}
\newtheorem{corollary}[theorem]{\bf Corollary}

\newtheorem{question}[theorem]{\bf Question}

\newtheorem{preproof}{{\bf Proof.}}

\begin{document}
\title{{\Large Ramsey numbers of 3-uniform loose paths and loose cycles}}
\author{\small  G.R. Omidi$^{\textrm{a},\textrm{b},1}$, M. Shahsiah$^{\textrm{a}}$\\
\small  $^{\textrm{a}}$Department of Mathematical Sciences,
Isfahan University
of Technology,\\ \small Isfahan, 84156-83111, Iran\\
\small  $^{\textrm{b}}$School of Mathematics, Institute for
Research
in Fundamental Sciences (IPM),\\
\small  P.O.Box:19395-5746, Tehran,
Iran\\
\small \texttt{E-mails: romidi@cc.iut.ac.ir,
m.shahsiah@math.iut.ac.ir}}
\date {}
\maketitle \footnotetext[1] {\tt This research was in part supported
by a grant from IPM (No.90050049)} \vspace*{-0.5cm}

\begin{abstract}
Haxell et. al. [%P. Haxell, T. Luczak, Y. Peng, V. R\"{o}dl, A.
%Ruci\'{n}ski, M. Simonovits, J. Skokan,
The Ramsey number for
hypergraph cycles I,  J. Combin.
 Theory, Ser. A, 113 (2006), 67-83] proved that the 2-color Ramsey number of
$3$-uniform loose cycles on $2n$ vertices is asymptotically
$\frac{5n}{2}$. Their proof is based on the method of Regularity
Lemma. Here, without using this method, we generalize their result
by determining the exact values of 2-color Ramsey numbers
involving loose paths and cycles in 3-uniform hypergraphs. More
precisely, we prove that for every $n\geq m\geq 3$,
$R(\mathcal{P}^3_n,\mathcal{P}^3_m)=R(\mathcal{P}^3_n,\mathcal{C}^3_m)=R(\mathcal{C}^3_n,\mathcal{C}^3_m)+1=2n+\lfloor\frac{m+1}{2}\rfloor$
and for $n>m\geq3$,
$R(\mathcal{P}^3_m,\mathcal{C}^3_n)=2n+\lfloor\frac{m-1}{2}\rfloor$.
These give a positive answer to a question of Gy\'{a}rf\'{a}s and
Raeisi [The Ramsey number of loose triangles and quadrangles in
hypergraphs,  Electron. J. Combin. 19 (2012), \#R30].\\

%Haxell et. al. in \cite{Ramsy number of loose cycle} proved that
%the 2-color Ramsey numbers of $3$-uniform loose cycles on $2n$
%vertices is asymptotically $\frac{5n}{2}$.
% Gy\'{a}rf\'{a}s et. al. \cite{subm} determined some off-diagonal
% 2-color Ramsey numbers involving loose paths and  cycles in 3-uniform
% hypergraphs when one of them has length at most six. Also they
%  conjectured that
% for every $n\geq m\geq 3$,
%$R(\mathcal{P}^3_n,\mathcal{P}^3_m)=R(\mathcal{P}^3_n,\mathcal{C}^3_m)=R(\mathcal{C}^3_n,\mathcal{C}^3_m)+1=2n+\Big\lfloor\frac{m+1}{2}\Big\rfloor.$
%In this paper we present a short proof of this conjecture. We also show that $R(\mathcal{P}^3_n,\mathcal{C}^3_m)=2m+\Big\lfloor\frac{n-1}{2}\Big\rfloor$,
%for $m>n\geq3$. \\\\

%Asymptotic values of 2-color Ramsey numbers for loose cycles and
%also loose paths were determined. Here we determine the exact
%value of  2-color Ramsey numbers involving loose paths and  cycles
%in 3-uniform hypergraphs, e. for every $n\geq m\geq 3$,
%$R(\mathcal{P}^3_n,\mathcal{P}^3_m)=R(\mathcal{P}^3_n,\mathcal{C}^3_m)=R(\mathcal{C}^3_n,\mathcal{C}^3_m)+1=2n+\Big\lfloor\frac{m+1}{2}\Big\rfloor.$
%This answers an open question of Gy\'{a}rf\'{a}s.


\end{abstract}

%%%%%%%%%%%%%%%%%%%%%%%%%%%%%%%%%%%%%%%%%%%%%%%%%%%%%%%%%%%%%%%%%%%%%%%%
\section{\normalsize Introduction}

\bigskip
A $k$-uniform {\it loose cycle} $\mathcal{C}_n^k$ (shortly, a {\it
cycle of length $n$})  is a hypergraph with vertex set
$\{v_1,v_2,\ldots,v_{n(k-1)}\}$ and with the set of $n$ edges
$e_i=\{v_1,v_2,\ldots, v_k\}+i(k-1)$, $i=0,1,\ldots, n-1$, where
we use mod $n(k-1)$ arithmetic and adding a number $t$ to a set
$H=\{v_1,v_2,\ldots, v_k\}$ means a shift, i.e. the set obtained
by adding $t$ to subscripts of each element of $H$. Similarly, a
$k$-uniform {\it loose path} $\mathcal{P}_n^k$ (simply, a {\it
path of length $n$}) is a hypergraph with vertex set
$\{v_1,v_2,\ldots,v_{n(k-1)+1}\}$  and with the set of $n$ edges
$e_i=\{v_1,v_2,\ldots, v_k\}+i(k-1)$, $i=0,1,\ldots, n-1$.
 For $k=2$ we get the usual definitions of a cycle $C_n$ and a path $P_n$ with $n$ edges. For an edge $e$ of a given loose path (also a given loose cycle)
$\mathcal{K}$, the first
 vertex and the last vertex are denoted by $v_{\mathcal{K},e}$ and
 $\hat{v}_{\mathcal{K},e}$, respectively. For two $k$-uniform hypergraphs $\mathcal{H}$ and $\mathcal{G}$,
  the \textit{Ramsey number} $R(\mathcal{H},\mathcal{G})$ is
  the smallest number $N$ such that each red-blue
coloring of the edges of the complete $k$-uniform hypergraph
$\mathcal{K}^k_N$ on $N$ vertices contains either a red copy of
$\mathcal{H}$ or a blue copy of $\mathcal{G}$.\\

One of the central problems in combinatorics and graph theory is
determining or estimating the Ramsey numbers. In this area, a
classical topic is the study of the Ramsey numbers of sparse
graphs,
 i.e., graphs with certain upper bound constraints on the degrees of the vertices.
 The study of these Ramsey numbers was initiated by Burr
and Erd\"{o}s \cite{Burr and erdos}. They
%asked for which graphs
%$H$ the
conjectured that for integer $\Delta$ there is a constant
$c=c(\Delta)$ so that for every graph $G$ on $n$ vertices and
maximum degree $\Delta$,  $R(G,G)\leq cn$.  This conjecture was
established  by Chv\'{a}tal et. al. in 1983 \cite{Ramsey number of
a graph with bounded maximum degree}. In this area the exact
values of $R(P_n,P_m)$, $R(P_n,C_m)$ and $R(C_n,C_m)$ are
determined (see \cite{Ramsey numbers for cycles,Ramsy number of
cycles,Ramsey number of cycle-path,Ramsey number of paths,cycle
ramsey numbers}).
 For a recent survey
including some results on the Ramsey numbers of
sparse graphs see \cite{survey}.\\

It was natural  to extend the results on the Ramsey number of
sparse graphs  for hypergraphs. In this direction, Kostochka and
R\"{o}dl \cite{Kostochka with given maximum degree} had shown that
for every $\epsilon>0$ and every positive integers $\Delta$ and
$r$, there exists $c=c(\epsilon,\Delta,r)$ such that
$R(\mathcal{H},\mathcal{H})$ for any $r$-uniform hypergraph
$\mathcal{H}$ with maximum degree at most $\Delta$ is at most
$c|V(\mathcal{H})|^{1+\epsilon}$; that is, the Ramsey numbers of
$k$-uniform hypergraphs of bounded maximum degree are almost
linear in their orders. In a sequel, Nagle et. al. \cite{Nagle
3-uniform hypergraphs} and Cooley et. al. \cite{Cooly 3-uniform
hypergraphs of bounded degree}  independently showed that the
Ramsey numbers of 3-uniform hypergraphs of bounded maximum degree
are linear in their order. This result was extended by several
authors to $k$-uniform hypergraphs.
%showed that
%for every $\epsilon >0$, the Ramsey number of any $k$-uniform
%hypergraph $H$ with $n$ vertices and maximum degree $\delta$
%satisfies
For an excellent article giving the latest results in this
direction, we refer the reader  to \cite{Ramsey numbers of sparse
hypergraphs}. Related with this topic,
  several interesting results
were obtained on 2-color Ramsey numbers of loose cycles and paths
in uniform hypergraphs. Haxell et. al. \cite{Ramsy number of loose
cycle} showed the following result on the Ramsey number of
3-uniform loose cycles, using Regularity Lemma.

\begin{theorem}\label{Haxell}

For all $\eta>0$ there exists $n_0=n_0(\eta)$ such that for every
$n> n_0$, every 2-colouring of $\mathcal{K}^{3}_{ 5(1+\eta)n/2}$
contains a monochromatic copy of $C_n^3$.
\end{theorem}

%Their proof is based on the Regularity Lemma.
Subsequently, Gy\'{a}rf\'{a}s,  S\'{a}rk\"{o}zy and Szemer\'{e}di
\cite{Ramsy number of loose cycle for k-uniform} extended this
result to $k$-uniform loose cycles and proved that for  $k\geq 3$
and all $\eta>0$ there exists $n_0=n_0(\eta)$ such that every
2-coloring of $\mathcal{K}^{k}_N$ with
$N=(1+\eta)\frac{1}{2}(2k-1)n$ contains a monochromatic copy of
$\mathcal{C}^k_n$, i.e.
 $R(\mathcal{C}^k_n,\mathcal{C}^k_n)$ is
asymptotically equal to $\frac{1}{2}(2k-1)n$. However, those
proofs are all based on application of the hypergraph regularity
method. Recently, Gy\'{a}rf\'{a}s and Raeisi determined the exact
values of
 2-color  Ramsey numbers of two k-uniform loose triangles and
 two k-uniform loose quadrangles \cite{subm}. They also posed the following
 question.

\begin{question}\label{question}
For every $n\geq m\geq 3$, is it true that
$R(\mathcal{P}^3_n,\mathcal{P}^3_m)=R(\mathcal{P}^3_n,\mathcal{C}^3_m)=R(\mathcal{C}^3_n,\mathcal{C}^3_m)+1=2n+\Big\lfloor\frac{m+1}{2}\Big\rfloor$?
In particular,
$$R(\mathcal{P}^3_n,\mathcal{P}^3_n)=R(\mathcal{C}^3_n,\mathcal{C}^3_n)+1=\Big\lceil\frac{5n}{2}\Big\rceil?$$
\end{question}
Related to question \ref{question}, it is proved \cite{The Ramsey
number of loose paths in 3-uniform hypergraphs} that for every
$n\geq \lfloor\frac{5m}{4}\rfloor$,
$R(\mathcal{P}^3_n,\mathcal{P}^3_m)=2n+\lfloor\frac{m+1}{2}\rfloor.$
In this article,  with a much short proof and without using the
Regularity Lemma, we give a positive answer to this question.
Moreover, we show that
$R(\mathcal{P}^3_m,\mathcal{C}^3_n)=2n+\lfloor\frac{m-1}{2}\rfloor$ for $n>m\geq 3$. Indeed, our results yield Theorem \ref{Haxell}. \\

\vspace{0.5cm}

The rest of this paper is organized as follows. In the next
section, we state the principal results required to prove the main
results. In section 3, we determine the exact value of the Ramsey
number of loose cycles in 3-uniform hypergraphs which is a
generalization of Theorem \ref{Haxell}. In section 4, we provide
the exact value of the Ramsey number of loose paths in 3-uniform
hypergraphs and finally in  the last section, the Ramsey number of
a loose path and a loose cycle in 3-uniform
hypergraphs is determined.\\\\
{\bf Note:} It is shown in \cite{subm} that
$(k-1)n+\lfloor\frac{m+1}{2}\rfloor$, $n\geq m \geq 2$ and $k\geq
3$, is a lower bound for $R(\mathcal{P}^k_n,\mathcal{P}^k_m)$,
$R(\mathcal{P}^k_n,\mathcal{C}^k_m)$ and
$R(\mathcal{C}^k_n,\mathcal{C}^k_m)+1$. Here we mention that  for
$n>m$ and $k\geq 3$, $R(\mathcal{P}^k_m,\mathcal{C}^k_n)\geq
(k-1)n+\lfloor\frac{m-1}{2}\rfloor$. To see that consider a
complete hypergraph whose vertex set is partitioned into two parts
$A$ and $B$, where $|A|=(k-1)n-1$ and
$|B|=\lfloor\frac{m-1}{2}\rfloor$. Color all edges that contain a
vertex of $B$ red, and the rest blue. In this coloring, the
longest red path has length at most $m-1$ and there is also no
blue copy of $\mathcal{C}_n^k$, since such a copy has $(k-1)n$
vertices. Our main aim in this article is to prove that these
lower bounds are tight. Therefore, in the rest of this paper, we
shall  prove just upper bounds for  the claimed Ramsey numbers.
Throughout the paper, we denote by $\mathcal{H}_{red}$ and
$\mathcal{H}_{blue}$ the induced 3-uniform hypergraph on edges of
color red and blue, respectively. Also we denote by
$|\mathcal{H}|$ and $\|\mathcal{H}\|$ the number of vertices and
edges of $\mathcal{H}$, respectively.

\section{\normalsize Preliminaries}

\medskip

In this section, we present some useful results.
% whose proofs are
%based on the property of monochromatic maximal paths in 3-uniform hypergraphs.
 We also rewrite some results from \cite{subm} in below.
%The first result that we need is a Lemma from \cite{subm}.

\medskip
\begin{lemma}\label{No Cn} \cite{subm}
Let $n\geq m\geq 3$ and
$\mathcal{K}^k_{(k-1)n+\lfloor\frac{m+1}{2}\rfloor}$ be 2-edge
colored red and blue. If $\mathcal{C}_{n}^k\subseteq
\mathcal{H}_{red}$, then either $\mathcal{P}_n^k\subseteq
\mathcal{H}_{red}$ or $\mathcal{P}_m^k\subseteq
\mathcal{H}_{blue}$. Also, if $\mathcal{C}_{n}^k\subseteq
\mathcal{H}_{red}$ then either $\mathcal{P}_n^k\subseteq
\mathcal{H}_{red}$ or $\mathcal{C}_m^k\subseteq
\mathcal{H}_{blue}$.
\end{lemma}
%\medskip
\begin{theorem}\label{R(Pk3,Pk3)} \cite{subm}
For every $k\geq 3$,\\\\
I)
$R(\mathcal{P}^k_3,\mathcal{P}^k_3)=R(\mathcal{C}^k_3,\mathcal{P}^k_3)=R(\mathcal{C}^k_3,\mathcal{C}^k_3)+1=3k-1$.\\\\
II)
$R(\mathcal{P}^k_4,\mathcal{P}^k_4)=R(\mathcal{C}^k_4,\mathcal{P}^k_4)=R(\mathcal{C}^k_4,\mathcal{C}^k_4)+1=4k-2$.

\end{theorem}
%%%%%%%%%%%%%%%%%%%%%%%%%%%%
%\medskip
%\begin{theorem}\label{R(Pk4,Pk4)} \cite{subm}
%For every $k\geq 3$,
%$R(\mathcal{P}^k_4,\mathcal{P}^k_4)=R(\mathcal{C}^k_4,\mathcal{P}^k_4)=R(\mathcal{C}^k_4,\mathcal{C}^k_4)+1=4k-2$.
%\end{theorem}
%%%%%%%%%%%%%%%%%%%%%%%%%%%%
\vspace{0.5 cm}

Let $\mathcal{P}$ be a loose path and $W$ be a set of  vertices
with $W\cap V(\mathcal{P})=\emptyset$. By a {\it
$\varpi_{\{v_i,v_j,v_k\}}$-configuration}, we mean a copy of
$\mathcal{P}^3_2$ with edges $\{x,v_i,v_j\}$ and \{$v_j,v_k,y\}$
so that $v_l$'s, $l\in\{i,j,k\}$, belong to two consecutive edges
of $\mathcal{P}$ and $\{x,y\}\subseteq W$. The vertices $x$ and
$y$ are called the end vertices of this configuration. A
$\varpi_{S}$-configuration, $S\subseteq V(e_i)\cup V(e_{i+1})$, is
good if   the last vertex of $e_{i+1}$ is not in $S$ and it is
bad, otherwise. Let $\mathcal{H}=\mathcal{K}_l^3$ be 2-edge
colored red and blue. We say that a red path
$\mathcal{P}=e_1e_2\ldots e_n$ of length $n$ is {\it maximal with
respect to} $W$ (for abbreviation w.r.t. $W$), for  $W\subseteq
V(\mathcal{H})\setminus V(\mathcal{P})$, if there are no vertices
$x$ and $y$ in $W$ so that for some $1\leq i\leq n$ either there
is a red path $\mathcal{P'}=e_1e_2\ldots e_{i-1}ee'e_{i+1}\ldots
e_n$ with $v_{\mathcal{P}',e}=v_{\mathcal{P},e_i}$ for $i=1$ and
$\hat{v}_{\mathcal{P}',e'}=\hat{v}_{\mathcal{P},e_i}$ for $i=n$,
or a red path $\mathcal{P'}=e_1e_2\ldots
e_{i-1}ee'e''e_{i+2}\ldots e_n$ with
$v_{\mathcal{P}',e}=v_{\mathcal{P},e_i}$ for $i=1$ and
$\hat{v}_{\mathcal{P}',e''}=\hat{v}_{\mathcal{P},e_{i+1}}$ for
$i=n-1$, in $\mathcal{H}$ so that
$V(\mathcal{P'})=V(\mathcal{P})\cup \{x,y\}$. We use these
definitions to deduce the following Lemma.

\bigskip
\begin{lemma}\label{spacial configuration}
Assume that $\mathcal{H}=\mathcal{K}^3_{n}$ is 2-edge colored red
and blue. Let $\mathcal{P}\subseteq \mathcal{H}_{red}$ be  maximal
w.r.t. $W\subseteq V(\mathcal{H})\setminus V(\mathcal{P})$ with
$|W|\geq 3$. Let  $e_i=\{v_{2i-1},v_{2i},v_{2i+1}\}$ and
$e_{i+1}=\{v_{2i+1},v_{2i+2},v_{2i+3}\}$ be consecutive edges of
$\mathcal{P}$. Either there is a good $\varpi_S$-configuration $C$
in $\mathcal{H}_{blue}$ with the set of end vertices in $W$ and
$S\subseteq V(e_i)\cup V(e_{i+1})$  or there is a bad
$\varpi_{S_1}$-configuration $C_1$, $S_1\subseteq V(e_i)\cup
V(e_{i+1})\setminus \{v_{2i+2}\}$. If there is no such a good
configuration $C$ in $\mathcal{H}_{blue}$ and
$e_{i+2}=\{v_{2i+3},v_{2i+4},v_{2i+5}\}$ is an edge of
$\mathcal{P}$, then also there is a  good
$\varpi_{S_2}$-configuration $C_2$, $S_2\subseteq V(e_{i+1})\cup
V(e_{i+2})$, in $\mathcal{H}_{blue}$ with end vertices in $W$
%such that $S_1\subseteq V(e_i)\cup V(e_{i+1})$,
 and $S_1\cap S_2=\emptyset$. Moreover, each vertex of $W$ except at
most one vertex can be considered as an end vertex of $C$ if there
exists such a configuration and otherwise, each vertex of $W$ can
be considered as an end vertex of $C_1$ and $C_2$.
\end{lemma}

%%%%%%%%%%%%%%%%%%%%%%%%%%%%
%%%%%%%%%%%%%%%%%%%%%%%%%%%%
\noindent\textbf{Proof. } Assume
 $W=\{x_1,...,x_t\}\subseteq V(\mathcal{K}^3_{n})\setminus V(\mathcal{P})$. If for some $x\in
W$, $\{v_{2i-1},v_{2i},x\}$ (resp.  $\{v_{2i+2},v_{2i+3},x\}$) is
red, then the maximality  of $\mathcal{P}$ w.r.t. $W$ implies that
for arbitrary vertices
 $x^{\prime}\neq x^{\prime\prime}\in W\setminus \{x\}$ the edges $\{x^{\prime},v_{2i+1},v_{2i}\}$ and
 $\{v_{2i},v_{2i+2},x^{\prime\prime}\}$ (resp. $\{x^{\prime},v_{2i+1},v_{2i+2}\}$ and $\{v_{2i+2},v_{2i},x^{\prime\prime}\}$)
 are blue and there is a good $\varpi_S$-configuration $C=\{x^{\prime},v_{2i+1},v_{2i}\}\{v_{2i},v_{2i+2},x^{\prime\prime}\}$
 (resp. $C=\{x^{\prime},v_{2i+1},v_{2i+2}\}\{v_{2i+2},v_{2i},x^{\prime\prime}\}$) with $S=\{v_{2i},v_{2i+1},v_{2i+2}\}\subseteq V(e_i)\cup V(e_{i+1})$.
 So we may assume that for every $x\in W$ both
  edges $\{v_{2i-1},v_{2i},x\}$ and $\{v_{2i+2},v_{2i+3},x\}$ are blue.
  If there is a vertex $y\in W$ such that at least one of the edges $f_1=\{v_{2i-1},v_{2i+1},y\}$,
  $f_2=\{v_{2i},v_{2i+1},y\}$,
  $f_3=\{v_{2i-1},v_{2i+2},y\}$ or $f_4=\{v_{2i},v_{2i+2},y\}$, say $f$,
  is blue, then there is a good $\varpi_S$-configuration $C=\{v_{2i-1},v_{2i},x\}f$, $x\neq y$, with
   $S=\{v_{2i-1},v_{2i}\}\cup f\setminus\{y\}\subseteq V(e_i)\cup V(e_{i+1})$.
   So if there is no such a good configuration $C$, then we may assume that for every $y\in
  W$ the edges  $f_1, f_2, f_3$ and $f_4$  are red.
 Therefore,
  maximality of $\mathcal{P}$ w.r.t. $W$ implies that for every $y'\in W$ the
  edge $\{v_{2i},v_{2i+3},y'\}$ is blue (otherwise,
replacing $e_ie_{i+1}$ by
$\{v_{2i-1},v_{2i+2},y\}\{y,v_{2i+1},v_{2i}\}\{v_{2i},v_{2i+3},y'\}$,
$y\neq y'$, in $\mathcal{P}$ yields a red path $\mathcal{P'}$
greater than $\mathcal{P}$, a contradiction). So, for  every
$a\neq b\in W$, $C_1=\{v_{2i-1},a,v_{2i}\}\{v_{2i},v_{2i+3},b\}$
is a bad $\varpi_{S_1}$-configuration with desired properties so
that $S_1=\{v_{2i-1},v_{2i},v_{2i+3}\}\subseteq V(e_i)\cup
V(e_{i+1})\setminus\{v_{2i+2}\}$.
%$S_1=\{v_1,v_2,v_5\}$ between $e_ie_{i+1}$ and $W$ with .
Clearly if there is $e_{i+2}=\{v_{2i+3},v_{2i+4},v_{2i+5}\}$ as an
edge of $\mathcal{P}$, then for every $x\in W$,
$\{v_{2i+2},v_{2i+4},x\}$ (also $\{v_{2i+1},v_{2i+4},x\}$) is
blue, otherwise replacing $e_ie_{i+1}$ by
$\{v_{2i-1},v_{2i+1},y\}\{y,v_{2i},v_{2i+2}\}\{v_{2i+2},v_{2i+4},x\}$
(replacing $e_ie_{i+1}$ by
$\{v_{2i-1},v_{2i+2},y\}\{y,$\\
$v_{2i},v_{2i+1}\}\{v_{2i+1},v_{2i+4},x\}$)
for some $y\neq x$ in $\mathcal{P}$ yields a red path
$\mathcal{P'}$ greater than $\mathcal{P}$, a contradiction to the
maximality of $\mathcal{P}$ w.r.t. $W$. So for every $a\neq b\in
W$, $C_2=\{a,v_{2i+2},v_{2i+4}\}\{v_{2i+4},v_{2i+1},b\}$ is a
 good $\varpi_{S_2}$-configuration with
desired properties where
$S_2=\{v_{2i+1},v_{2i+2},v_{2i+4}\}\subseteq V(e_{i+1})\cup
V(e_{i+2})$.
% with $S_2=\{v_3,v_4,v_6\}$ between $e_{i+1}e_{i+2}$ and $W$
 $\hfill\blacksquare$\\

%%%%%%%%%%%%%%%%%%%%%%%%%%%%%%%%%%%%%%%%%%%%%%%%%%%%%%%%%%%%%%%%%%%%%%%%%%%%%%%%%%%%%%%%%%%%%%%%%%%%%%
%Using the same argument to the proof of Lemma \ref{spacial
%configuration} we have the following result.
%The followings are  immediate corollaries of the  Lemma
%\ref{spacial configuration}.
% and guarantee the existence of

%%%%%%%%%%%%%%%%%%%%%%%%%%%%

 The following is an  immediate corollary
of the  Lemma \ref{spacial configuration}.

\begin{corollary}\label{there is a P2 or P4}
Assume that $\mathcal{H}=\mathcal{K}^3_{n}$ is 2-edge colored red
and blue. Let $\mathcal{P}\subseteq \mathcal{H}_{red}$ be  maximal
w.r.t.  $W\subseteq V(\mathcal{H})\setminus V(\mathcal{P})$ with
$|W|\geq 3$ . Let $e_i=\{v_{2i-1},v_{2i},v_{2i+1}\}$,
$e_{i+1}=\{v_{2i+1},v_{2i+2},v_{2i+3}\}$ and
$e_{i+2}=\{v_{2i+3},v_{2i+4},v_{2i+5}\}$ be consecutive edges of
$\mathcal{P}$. Either for every distinct vertices $x$ and $y$ of
$W$, except at most one vertex there is a blue path
$Q=\{\bar{v}_1,x,\bar{v}_2\}\{\bar{v}_2,y,\bar{v}_3\}$ of length 2
with $\{\bar{v}_1,\bar{v}_2,\bar{v}_3\}\subseteq V(e_{i})\cup
V(e_{i+1})\setminus \{v_{2i+3}\}$ or for every distinct vertices
$x$,$y$ and $z$ of $W$ there is a blue path
$Q'=\{v'_1,x,v'_2\}\{v'_2,v'_3,y\}\{y,v'_4,v'_5\}\{v'_5,v'_6,z\}$
of length 4 with $\{v'_1,v'_2,v'_3,v'_4,v'_5,v'_6\}\subseteq
V(e_{i})\cup V(e_{i+1})\cup V(e_{i+2})\setminus \{v_{2i+5}\}$.
%so that each vertex of $W$ except at most
%one vertex can be considered as an end vertex of $Q$ and each
%vertex of $W$ can be considered as an end vertex of $Q'$.
\end{corollary}

\begin{corollary}\label{there is a Pl}
Let $\mathcal{H}=\mathcal{K}_l^3$ be two edge colored red and blue
and $\mathcal{P}=e_1e_2\ldots e_n$ be maximal red path w.r.t. $W$,
$W\subseteq V(\mathcal{H})\setminus V(\mathcal{P})$. Then for some
$r\geq 0$ and $W'\subseteq W$ there is a blue path $Q$ between
$W'$ and $\bar{\mathcal{P}}=e_1e_2\ldots e_{n-r}$ with no
$\hat{v}_{\mathcal{P}, e_{n-r}}$ as a vertex and
$\|Q\|=2(|W'|-1)\geq n-r$ where either $x=|W\setminus W'|=0,1$ or
$x\geq 2$ and $0\leq r \leq 2$.
%for every $W\subseteq V(\mathcal{H})\setminus V(\mathcal{P})$.
%Then there is a blue path $Q$ of length
%$t=min\{2\lfloor\frac{n}{2} \rfloor, 2|W|-??\}$ between
%$e_1e_2\ldots e_t$ and $W$ so that the end vertices of $Q$ are in
%$W$ and each vertex of $W$, except end vertices, lies on two edges
%of $Q$. Moreover the end vertices of $Q$ ???.
\end{corollary}

%%%%%%%%%%%%%%%%%%%%%%%%%%%%
\noindent\textbf{Proof. }Let $\mathcal{P}=e_1e_2\ldots e_{n}$ be a
maximal  red path w.r.t. $W=V(\mathcal{H})\setminus
V(\mathcal{P})$ where $e_i=\{v_1,v_2,v_3\}+2(i-1)$, $i=1,\ldots,
n$. Since $\mathcal{P}_1=\mathcal{P}$ is maximal w.r.t. $W_1=W$,
using Corollary
 \ref{there is a P2 or P4} there is a blue path $Q_1$ with end vertices $x_1$ and $y_1$ in $W_1$ between $E'_1$ and $W_1$
 where $E_1=\{e_i: i_1=1\leq i \leq 2\}$, $\bar{E}_1=E_1\cup\{e_{i_1+2}\}$, $E'_1\in \{E_1,\bar{E}_1\}$, $\| Q_1\|=2|E'_1|-2$
 and $Q_1$ does not contain the last vertex of $E'_1$. Set $i_2=min\{j:
j\in\{i_1+2,i_1+3\},
 e_j\not\in E'_1\}$, $E_2=\{e_i: i_2\leq i\leq i_2+1\}$,
$\bar{E}_2=E_2\cup
 \{e_{i_2+2}\}$, $\mathcal{P}_2=\mathcal{P}_1\setminus E'_1$ and $W_2=(W\setminus V(Q_1))\cup \{x_1,y_1\}$.
 Since $\mathcal{P}_2$ is maximal w.r.t. $W_2$,
 using Corollary \ref{there is a P2 or P4} there is a blue path $Q_2$ between
 $E'_2$, $E'_2\in \{E_2, \bar{E}_2\}$, and $W_2$ such that $\| Q_2\|=2|E'_2|-2$
 and $Q_2$ does not contain the last vertex of $E'_2$ and
  $Q_1\cup Q_2$ is a blue path with end
 vertices $x_2, y_2$ in $W_2$.
  Now set $i_3=min\{j: j\in\{i_2+2,i_2+3\}, e_j\not\in E'_2\}$, $E_3=\{e_i: i_3\leq i\leq
 i_3+1\}$, $\bar{E}_3=E_3\cup
 \{e_{i_3+2}\}$, $\mathcal{P}_3=\mathcal{P}_2\setminus E'_2$ and $W_3=(W\setminus V(Q_1\cup Q_2))\cup
 \{x_2,y_2\}$ and continue this process. Assume that this process is terminated in $t-$th
 step. Clearly either $|W\setminus W_t|=0,1$ or $|W\setminus W_t|\geq
 2$ and $0\leq \|\mathcal{P}_t\|\leq 2$. So $Q=Q_1\cup Q_2\cup \ldots \cup
 Q_{t-1}$ is a blue path between $\bar{\mathcal{P}}=e_1e_2\ldots
 e_{n-r}$ and $W'=(W\setminus W_t)\cup\{x_{t-1},y_{t-1}\}$ where $\mathcal{P}_t=e_{n-r+1}e_{n-r+2}\ldots
 e_{n}$.
 $\hfill\blacksquare$
%%%%%%%%%%%%%%%%%%%%%%%%%%%%%%%%%%%%%%%%%%%%%%%%%%%%%%%%%%%%%%%%%%%%%%%%%%%%%%%%%%%%%%%%%%%%%%%%%%%%%%
\section{\normalsize Cycle-cycle Ramsey number in 3-uniform hypergraphs}

%%%%%%%%%%%%%%%%%%%%%%%%%%%%%%%%%%%%%%%%%%%%%%%%%%%%%%%%%%%%%%%%%%%%%%%%%%%%%%%%%%%%%%%%%%%%%%%%%%%%%%
In this section, we provide the exact value of
$R(\mathcal{C}^3_n,\mathcal{C}^3_m)$, when  $n\geq m\geq 3$.
Before that we need two following lemmas.

 \bigskip
\begin{lemma}\label{cn-1 implies cm}
Let $n\geq m\geq 3$, $(n,m)\neq (3,3),(4,3),(4,4)$ and
$\mathcal{K}^3_{2n+\lfloor\frac{m-1}{2}\rfloor}$ be 2-edge colored
red and blue. Assume there is no copy of $\mathcal{C}^3_{n}$ in
$\mathcal{H}_{red}$ and $\mathcal{C}=\mathcal{C}^3_{n-1}$ is a
loose cycle in $\mathcal{H}_{red}$. Then there is a copy of
$\mathcal{C}^3_{m}$ in $\mathcal{H}_{blue}$.  Moreover, for $n>m$,
there is also a copy of $\mathcal{P}^3_m$ in $\mathcal{H}_{blue}$.
\end{lemma}
%%%%%%%%%%%%%%%%%%%%%%%%%%%%
\noindent\textbf{Proof. }Let $t=2n+\lfloor\frac{m-1}{2}\rfloor$
and $\mathcal{C}=e_1e_2\ldots e_{n-1}$ be a copy of
$\mathcal{C}_{n-1}^3$ in $\mathcal{H}_{red}$ with  edges
$e_i=\{v_1,v_2,v_3\}+2(i-1)$ (mod $2(n-1)$), $i=1,\ldots, n-1$.
Let $W=V(\mathcal{H})\setminus V(\mathcal{C})$ where
$\mathcal{H}=\mathcal{K}^3_t$.

\bigskip
\noindent \textbf{Case 1. } There are an edge
$e_i=\{v_{2i-1},v_{2i},v_{2i+1}\}$, $1\leq i\leq n-1$, and a
vertex $z\in W$ such that $\{v_{2i},v_{2i+1},z\}$ is red.

\medskip
Let $\mathcal{P}=e_{i+1}e_{i+2}\ldots e_{n-1}e_1e_2\ldots
e_{i-2}e_{i-1}$ and $W_0=W\setminus \{z\}$. Since $\mathcal{P}$ is
a
 maximal path w.r.t.  $W_0$, using Corollary \ref{there is a
 Pl}, there is a  blue path $Q$ of length $l'$ between
 $\bar{\mathcal{P}}$, the path obtained from $\mathcal{P}$ by deleting the last $r$ edges, and $W'$ for
 some $r\geq 0$ and $W'\subseteq W_0$ with mentioned properties in Corollary \ref{there is a Pl}.
%   the set of $\{\bar{c}_1,\bar{c}_2,\ldots,\bar{c}_{q'}\}$ of
%configurations.
Let $x',y'$  be the end vertices of $Q$  in $W'$,
 $T=W_0\setminus V(Q)$ and $x=|T|$. We have following subcases.

\medskip
{\it Subcase 1}. $x=0.$

\medskip
Clearly  $l'=2\lfloor\frac{m-1}{2}\rfloor$. First let
 $m$ be even. Hence  $l'=m-2$ and so
$Q\{y',v_{2i},v_{2i-1}\}$\\
$\{v_{2i-1},z,x'\}$ is a blue $\mathcal{C}^3_m$. Moreover, if
$n>m$, then $r\geq 1$ and so at least one
 of $\{v_{2i-3},v_{2i-2},x'\}Q\{y',v_{2i-1},v_{2i}\}$ or
$Q\{y',v_{2i-1},v_{2i}\}\{v_{2i},v_{2i-2},z\}$ is a blue copy of
$\mathcal{P}_{m}^3$. Now let $m$ be odd. So $l'=m-1$. In this case
we truncate $Q$ to a path $Q'$ of length $m-3$ so that
$v_{2i-2}\notin Q'$ by removing the last two edges. Now we may
assume the vertices $x'$ and $y''\neq y'$ of $W_0$ are the end
vertices of $Q'$.
 So $Q'\{y'',v_{2i-2},v_{2i}\}\{v_{2i},y',v_{2i-1}\}\{v_{2i-1},z,x'\}$
is a copy of $\mathcal{C}^3_m$ in $\mathcal{H}_{blue}$.  Also
$Q\{y',v_{2i-1},v_{2i}\}$ is a blue copy of  $\mathcal{P}_{m}^3$.

\medskip
{\it Subcase 2}. $x=1.$

\medskip
Let $T=\{u\}$. Clearly $l'=2\lfloor\frac{m-1}{2}\rfloor-2$. Let
$m$ be odd. Then $l'=m-3$ and $r\geq 1$. So
$Q\{y',v_{2i-2},v_{2i}\}\{v_{2i},u,v_{2i-1}\}\{v_{2i-1},z,x'\}$ is
a blue copy of $\mathcal{C}^3_m$. If $n>m$, then $r\geq 2$ and
since $\hat{\mathcal{P}}=e_{i-2}e_{i-1}e_{i}$ is a maximal path
w.r.t. $\hat{W}=\{x',y',u,z\}$, using Corollary \ref{there is a P2
or P4} there is either a blue path $Q'$ of length 2 between
$(V(e_{i-2})\cup V(e_{i-1}))\setminus \{v_{2i-1}\}$ and $\hat{W}$
or a blue path $Q'$ of length 4 between
$V(\hat{\mathcal{P}})\setminus \{v_{2i+1}\}$ and $\hat{W}$ so that
$Q\cup Q'$ is a blue path, say $Q''$.  Now let $l''$ be the length
of $Q'$. If $l''=4$, $Q''$ is a blue $\mathcal{P}_{m+1}^3$ and so
there is a $\mathcal{P}_{m}^3$ in $\mathcal{H}_{blue}$. If
$l''=2$, the length of $Q''$ is $m-1$. Without lose of generality
let $x',y''$ be the end vertices of $Q''$. Then
$\{v_{2i-1},v_{2i},x'\}Q''$ is a blue copy of $\mathcal{P}_{m}^3$.
 Now suppose $m$ is even, so
$l'=m-4$ and $r\geq 2$. If
  $\{v_{2i-5},v_{2i-4},x'\}$ is red, then $Q\{y',v_{2i-3},v_{2i-4}\}\{v_{2i-4},v_{2i-2},u\}\{u,v_{2i},v_{2i-1}\}\{v_{2i-1},z,x'\}$
is a blue $\mathcal{C}_m^3$. If $\{v_{2i-4},v_{2i-3},u\}$ is red,
then
$Q\{y',v_{2i-1},v_{2i}\}\{v_{2i},v_{2i-2},u\}\{u,z,v_{2i-5}\}\{v_{2i-5},v_{2i-4},$\\
$x'\}$ is a  blue copy of $\mathcal{C}_m^3$. Otherwise,
$Q\{y',v_{2i-1},v_{2i}\}\{v_{2i},v_{2i-2},u\}\{u,v_{2i-3},v_{2i-4}\}$\\
$\{v_{2i-4},v_{2i-5},x'\}$ is a
 copy of $\mathcal{C}_m^3$ in $\mathcal{H}_{blue}$.
For $n>m$, clearly $r\geq 3$ and since
$\hat{\mathcal{P}}=e_{i-3}e_{i-2}e_{i-1}$ is a maximal path w.r.t.
$\hat{W}=\{x',y',u,z\}$, using Corollary \ref{there is a P2 or P4}
there is a blue path $Q'$ of length either $l''=2$, between
$(V(e_{i-3})\cup V(e_{i-2}))\setminus \{v_{2i-3}\}$ and $\hat{W}$
or $l''=4$, between $V(\hat{\mathcal{P}})\setminus \{v_{2i-1}\}$
and $\hat{W}$ such that $Q\cup Q'$ is a blue path, say $Q''$. If
$l''=4$, then $Q''$ is a blue $\mathcal{P}_m^3$. otherwise the
length of $Q''$ is $m-2$ and we may assume that $x',y''$ be the
end vertices of $Q''$. If $y''=z$,
 then $\{u,v_{2i-2},v_{2i}\}\{v_{2i},v_{2i-1},x'\}Q''$ is a blue copy of  $\mathcal{P}_{m}^3$.
 Otherwise, at least one
 of $\{v_{2i-3},v_{2i-2},x'\}Q''\{y'',v_{2i-1},v_{2i}\}$ or
$Q''\{y'',v_{2i-1},v_{2i}\}\{v_{2i},v_{2i-2},z\}$ is a blue copy
of $\mathcal{P}_{m}^3$.\\

\medskip
{\it Subcase 3}. $x\geq 2$.

\medskip
One can easily check that $r\geq 3$. This case does not occur by
Corollary \ref{there is a Pl}.

\bigskip
\noindent \textbf{Case 2. } There are
$e_i=\{v_{2i-1},v_{2i},v_{2i+1}\}$, $1\leq i\leq n-1$, and a
vertex $z\in W$ such that $\{v_{2i-1},v_{2i},z\}$ is red.

\medskip
In this case, consider the path $ \mathcal{P}=e_{i-1}e_{i-2}\ldots
e_2e_1e_{n-1}e_{n-2}\ldots e_{i+2}e_{i+1}$ and repeat the proof of
case 1. By
  an argument similar to the case 1 we can find a blue copy of $\mathcal{C}^3_m$ and a blue copy of $\mathcal{P}^3_m$ for $n>m$.

\bigskip
\noindent \textbf{Case 3. } For every
$e_i=\{v_{2i-1},v_{2i},v_{2i+1}\}$, $1\leq i\leq n-1$, and every
vertex $z\in W$  the edges $\{v_{2i-1},v_{2i},z\}$ and
$\{v_{2i},v_{2i+1},z\}$ are blue.

\medskip
 In this case, clearly there are
blue copies of  $\mathcal{C}^3_m$ and $\mathcal{P}^3_m$ between
$V(\mathcal{C})$ and $W$ and these observations complete the
proof.
 $\hfill\blacksquare$

  \bigskip
\begin{lemma}\label{R(C43,C33)}
$R(\mathcal{C}^3_4,\mathcal{C}^3_3)=9$
\end{lemma}

%%%%%%%%%%%%%%%%%%%%%%%%%%%%
\noindent\textbf{Proof. } Let $\mathcal{H}=\mathcal{K}^3_9$ be
2-edge colored red and blue. Suppose that there is no red copy of
$\mathcal{C}^3_4$ and no blue copy of $\mathcal{C}^3_3$. Using
Theorem \ref{R(Pk3,Pk3)} we may assume that there is a blue copy
of $\mathcal{C}^3_4$. Let $\mathcal{C}=e_1e_2e_3e_4$ be a copy of
$\mathcal{C}_4^3$ in $\mathcal{H}_{blue}$ with  edges
$e_i=\{v_1,v_2,v_3\}+2(i-1)$ (mod $8$), $i=1,\ldots, 4$. Let $v\in
V(\mathcal{H})\setminus V(\mathcal{C})$.
 Since
there is no blue copy of $\mathcal{C}^3_3$,
$\{v_1,v_6,v\}\{v,v_3,v_8\}\{v_8,v_4,v_2\}\{v_2,v_5,v_1\}$ is a
red copy of $\mathcal{C}^3_4$, a contradiction.
 $\hfill\blacksquare$\\

The main result of this section is the following  result on the
Ramsey number of loose cycles in 3-uniform hypergraphs.

\bigskip
\begin{theorem}\label{main theorem3}
For every $n\geq m\geq 3$,
$$R(\mathcal{C}^3_n,\mathcal{C}^3_m)=2n+\Big\lfloor\frac{m-1}{2}\Big\rfloor.$$
\end{theorem}
%%%%%%%%%%%%%%%%%%%%%%%%%%%%
%%%%%%%%%%%%%%%%%%%%%%%%%%%%
\noindent\textbf{Proof. } We prove this theorem by induction on
$m+n$.   Using Theorem \ref{R(Pk3,Pk3)}, the base case is trivial.
%Suppose that for $m^{\prime}+n^{\prime}< m+n$ with $n'\geq m'\geq
%3$, $R(\mathcal{C}^3_{n^{\prime}},
%\mathcal{C}^3_{m^{\prime}})=2n^{\prime}+\Big\lfloor\frac{m^{\prime}-1}{2}\Big\rfloor.$
%Now, let $n\geq m\geq 3$.
 Suppose indirectly that the edges of $\mathcal{H}=\mathcal{K}^3_{2n+\lfloor\frac{m-1}{2}\rfloor}$ can be
 colored red and blue without a red copy of $\mathcal{C}^3_n$ and
 a blue copy of $\mathcal{C}^3_m$.
 Consider the following cases.

\bigskip
\noindent \textbf{Case 1. } $n=m$

\medskip
Using Theorem \ref{R(Pk3,Pk3)} and Lemma \ref{R(C43,C33)} we may
assume that $n\geq 5$. By the induction hypothesis,
$R(\mathcal{C}^3_{n-1},\mathcal{C}^3_{n-1})=
2(n-1)+\Big\lfloor\frac{n-2}{2}\Big\rfloor<
2n+\Big\lfloor\frac{n-1}{2}\Big\rfloor$. So we may assume that
there is a red copy of $\mathcal{C}^3_{n-1}$ in $\mathcal{H}$.
Using Lemma \ref{cn-1 implies cm} we have a blue
$\mathcal{C}^3_{n}$, a contradiction.

\bigskip
\noindent \textbf{Case 2. }$n> m$.

\medskip
 In this case, $n-1\geq m$ and since
 $R(\mathcal{C}^3_{n-1},\mathcal{C}^3_{m})=
2(n-1)+\Big\lfloor\frac{m-1}{2}\Big\rfloor<
2n+\Big\lfloor\frac{m-1}{2}\Big\rfloor$, we may assume that
$\mathcal{C}^3_{n-1}\subseteq \mathcal{H}_{red}$.
 Using  Lemmas \ref{cn-1 implies cm} we have a blue
$\mathcal{C}^3_{m}$, a contradiction. These observations complete
the proof.
 $\hfill\blacksquare$
%%%%%%%%%%%%%%%%%%%%%%%%%%%%%%%%%%%%%%%%%%%%%%%%%%%%%%%%%%%%%%%%%%%%%%%%%%%%%%%%%%%%%%%%%%%%%%%%%%%%%%%%%%%

\section{\normalsize Path-path Ramsey number in 3-uniform hypergraphs}

\bigskip
In this section, we determine the exact value of
$R(\mathcal{P}^3_n,\mathcal{P}^3_m)$, for  $n\geq m\geq 3$. To see
that we need the following Lemmas.

\bigskip
\begin{lemma}\label{pn-1 implies cm and pm}
Let $n\geq m\geq 3$, $(n,m)\neq (3,3),(4,3),(4,4)$ and
$\mathcal{K}^3_{2n+\lfloor\frac{m+1}{2}\rfloor}$ be 2-edge colored
red and blue. If $\mathcal{P}=\mathcal{P}^3_{n-1}$ is a maximum
red path, then there is a copy of $\mathcal{P}^3_{m}$ in
$\mathcal{H}_{blue}$.
\end{lemma}
%%%%%%%%%%%%%%%%%%%%%%%%%%%%
\noindent\textbf{Proof. }Let $t=2n+\lfloor\frac{m+1}{2}\rfloor$,
$\mathcal{P}=e_1e_2\ldots e_{n-1}$ be a red  copy of
$\mathcal{P}_{n-1}^3$ with  edges $e_i=\{v_1,v_2,v_3\}+2(i-1)$,
$i=1,\ldots, n-1$ and $\bar{W}=V(\mathcal{K}^3_t)\setminus
 V(\mathcal{P})$.  By Lemma \ref{No Cn}, we
 may assume that there is no  copy of
 $\mathcal{C}^3_{n}$ in $\mathcal{H}_{red}$. Let $W_0=\bar{W}\setminus \{u\}$ for $u\in \bar{W}$. Since $\mathcal{P}'=\mathcal{P}\setminus \{e_1\}$ is a
 maximal path w.r.t. $W_0$, using  Corollary \ref{there is a
 Pl}, there is a
%  and $n-2=2l+r$ where $0\leq r\leq
% 1$.  Partition the set $E(\mathcal{P})\setminus \{e_1\}$ into $l$ classes $A_1,A_2,\ldots,A_l$ such
%  that each class contains two consecutive edges of $\mathcal{P}$. Using Lemma \ref{spacial configuration}, we can find a
  blue path $Q=\mathcal{P}_{l'}$, with no $v_{2n-1}$ as a vertex,  of length $l'=2q'$ between
 $\bar{\mathcal{P}}'=e_2e_3\ldots e_{n-1-r}$ and $W'$ for some $r\geq 0$ and $W'\subseteq W_0$ with the mentioned properties in Corollary  \ref{there is a Pl}.
%   the set of $\{\bar{c}_1,\bar{c}_2,\ldots,\bar{c}_{q'}\}$ of
%configurations.
Let $y,z$  be the end vertices of $Q$  in $W'$, $T=W_0\setminus
W'$ and $x=|T|$.
% Clearly $|T|=\lfloor\frac{m+1}{2}\rfloor-q'$. Assume $m=2k+r'$ for some
%$r'$, $0\leq r'\leq 1$.
 We have one of the following cases.

\bigskip
\noindent \textbf{Case 1. }$x=0$.

\medskip
 It is easy to see that
$l'=2\lfloor\frac{m+1}{2}\rfloor-2$. If $m$ is odd, $l'=m-1$ and
clearly $Q\{y,v_1,u\}$ is a blue copy of $\mathcal{P}_m^3$. If $m$
is even, $l'=m-2$. Since $v_{2n-1}$ is not a vertex of $Q$ and
there is no copy of red $\mathcal{C}_n^3$, clearly $\{v_1,u,y\}Q
\{z,v_2,v_{2n-1}\}$ is a blue copy of $\mathcal{P}_m^3$.

\bigskip
\noindent \textbf{Case 2. }$x=1$.

\medskip
Let $T=\{v\}$. In this case, $l'=2\lfloor\frac{m+1}{2}\rfloor-4$.
Clearly for odd $m$, $l'=m-3$ and $r\geq 1$. One can easily check
that $Q\{z,v_2,v_{2n-2}\}\{v_{2n-2},u,v\}\{v,v_{2n-1},v_1\}$ is a
blue $\mathcal{P}_m^3$. For even $m$, clearly $l'=m-4$ and $r\geq
2$. Since $\hat{\mathcal{P}}=e_{n-2}e_{n-1}$ is maximal w.r.t.
$\hat{W}=\{y,z,u,v\}$, by using  Lemma \ref{spacial
configuration}, there is a blue $\varpi_S$-configuration, say
$\mathcal{P}''$, with $S\subseteq V(e_{n-2})\cup V(e_{n-1})$ so
that $Q'=Q\cup
 \mathcal{P}''$ is a blue path of length $m-2$ and at least one of
 $v_{2n-2}$ and $v_{2n-1}$, say $w$, is not in $V(Q')$. Without lose of
 generality assume that y and $v$ are the end vertices of $Q'$. So
 $\{u,v_1,y\}Q'\{v,w,v_2\}$ is a blue
 $\mathcal{P}_m^3$.

\bigskip
\noindent \textbf{Case 3. }$x\geq 2$.

\medskip
One can easily check that $r\geq 3$. This case does not  occur by
Corollary \ref{there is a Pl}.
 $\hfill\blacksquare$

\vspace{0.5 cm} Using Theorem \ref{R(Pk3,Pk3)} we have
$R(\mathcal{P}^3_4,\mathcal{P}^3_4)=10$. Since
$R(\mathcal{P}^3_4,\mathcal{P}^3_3)\leq
R(\mathcal{P}^3_4,\mathcal{P}^3_4)$, we have the following Lemma.

\begin{lemma}\label{R(P34,P33)}
$R(\mathcal{P}^3_4,\mathcal{P}^3_3)=10.$
\end{lemma}
%%%%%%%%%%%%%%%%%%%%%%%%%%%%%%%%%%%%%%%%%%%%%%%%%%%%%%%%%%%%%%%%%%%%%%%%%%%%%%%%%%%%%%%%%%%%%%%%%%%%%%

%%%%%%%%%%%%%%%%%%%%%%%%%%%%%%%%%%%%%%%%%%%%%%%%%%%%%%%%%%%%%%%%%%%%%%%%%%%%%%%%%%%%%%%%%%%%%%%%%%%%%%

%%%%%%%%%%%%%%%%%%%%%%%%%%%%%%%%%%%%%%%%%%%%%%%%%%%%%%%%%%%%%%%%%%%%%%%%%%%%%%%%%%%%%%%%%%%%%%%%%%%%%%
%%%%%%%%%%%%%%%%%%%%%%%%%%%%%%%%%%%%%%%%%%%%%%%%%%%%%%%%%
%%%%%%%%%%%%%%%%%%%%%%%%%%%%%%%%%%%%%%%%%%%%%%%%%%%%%%%%%%%%%%%%%%%%%%%%%%%%%%%%%%%%%%%%%%%%%%%%%%%%%%
\bigskip
\begin{theorem}\label{main theorem1}
For every $n\geq m\geq 3$,
$$R(\mathcal{P}^3_n,\mathcal{P}^3_m)=2n+\Big\lfloor\frac{m+1}{2}\Big\rfloor.$$
\end{theorem}
%%%%%%%%%%%%%%%%%%%%%%%%%%%%
%%%%%%%%%%%%%%%%%%%%%%%%%%%%
\noindent\textbf{Proof. } We give a proof  by induction on $m+n$.
Using Theorem \ref{R(Pk3,Pk3)} the base case is trivial.
% Suppose that for $m^{\prime}+n^{\prime}< m+n$ with $n'\geq m'\geq 3$,
%$R(\mathcal{P}^3_{n^{\prime}},
%\mathcal{P}^3_{m^{\prime}})=2n^{\prime}+\Big\lfloor\frac{m^{\prime}+1}{2}\Big\rfloor$.
%Now, let $n\geq m\geq 3$ and let
Let $\mathcal{H}=\mathcal{K}^3_{2n+\lfloor\frac{m+1}{2}\rfloor}$
be 2-edge colored red and blue with  no a red copy
$\mathcal{P}^3_n$ of and a  blue copy of $\mathcal{P}^3_m$.
Consider the following cases.

\bigskip
\noindent \textbf{Case 1. } $n=m$

\medskip

By Theorem \ref{R(Pk3,Pk3)} and Lemma \ref{R(P34,P33)}  we may
assume $n\geq 5$. By the induction hypothesis,
$R(\mathcal{P}^3_{n-1},\mathcal{P}^3_{n-1})=
2(n-1)+\Big\lfloor\frac{n}{2}\Big\rfloor<
2n+\Big\lfloor\frac{n+1}{2}\Big\rfloor$. So we may assume that
there is a red copy of $\mathcal{P}^3_{n-1}$. Using Lemma
\ref{pn-1 implies cm and pm} we have a blue copy of
$\mathcal{P}^3_{n}$ in $\mathcal{H}$, a contradiction.

\bigskip
\noindent \textbf{Case 2. }$n> m$.

\medskip
 In this case, $n-1\geq m$. Since
 $R(\mathcal{P}^3_{n-1},\mathcal{P}^3_{m})=
2(n-1)+\Big\lfloor\frac{m+1}{2}\Big\rfloor<
2n+\Big\lfloor\frac{m+1}{2}\Big\rfloor$, we may assume there is a
copy of $\mathcal{P}^3_{n-1}$ in $\mathcal{H}_{red}$.
 Using  Lemma \ref{pn-1 implies cm and pm} we have a blue copy of
$\mathcal{P}^3_{m}$ in $\mathcal{H}$,a contradiction. These
observations complete the proof.
 $\hfill\blacksquare$

%%%%%%%%%%%%%%%%%%%%%%%%%%%%%%%%%%%%%%%%%%%%%%%%%%%%%%%%%%%%%%%%%%%%%%%%%%%%%%%%%%%%%%%%%%%%%%%%%%%%%%%%%%%

\section{\normalsize Path-cycle Ramsey numbers in 3-uniform hypergraphs }
In this section, the Ramsey number of  a loose path and a loose
cycle in 3-uniform hypergraphs is determined.

 \vspace{0.5 cm} It is worth noting that  we can
also conclude $R(\mathcal{P}^3_n,\mathcal{C}^3_m)\leq
n+\lfloor\frac{m+1}{2}\rfloor$, for $n \geq m\geq 3$. To see that
assume
$\mathcal{H}=\mathcal{K}^3_{2n+\lfloor\frac{m+1}{2}\rfloor}$ be
2-edge colored red and blue with  no a red copy of
$\mathcal{P}^3_n$ and no a  blue copy of $\mathcal{C}^3_m$. Since,
using Theorem \ref{main theorem3},
$R(\mathcal{C}^3_n,\mathcal{C}^3_m)=2n+\lfloor\frac{m-1}{2}\rfloor<
2n+\lfloor\frac{m+1}{2}\rfloor$, we have a red copy of
$\mathcal{C}^3_n$ in $\mathcal{H}$. Therefore, Lemma \ref{No Cn}
implies a contradiction to our assumptions. So we have the
following theorem.

 \bigskip
\begin{theorem}\label{main theorem2}
For every $n\geq m\geq 3$,
$$R(\mathcal{P}^3_n,\mathcal{C}^3_m)=2n+\Big\lfloor\frac{m+1}{2}\Big\rfloor.$$
\end{theorem}
%%%%%%%%%%%%%%%%%%%%%%%%%%%%
%%%%%%%%%%%%%%%%%%%%%%%%%%%%
%\noindent\textbf{Proof.} Suppose to the contrary that
%$\mathcal{H}=\mathcal{K}^3_{2n+\lfloor\frac{m+1}{2}\rfloor}$ is
%2-edge colored red and blue with  no a red copy of
%$\mathcal{P}^3_n$  and a blue copy of $\mathcal{C}^3_m$. Since,
%using Theorem \ref{main theorem3},
%$R(\mathcal{C}^3_n,\mathcal{C}^3_m)=2n+\lfloor\frac{m-1}{2}\rfloor<
%2n+\lfloor\frac{m+1}{2}\rfloor$, we have a red copy of
%$\mathcal{C}^3_n$. Therefore, Lemma \ref{No Cn} implies a
%contradiction to our assumptions. This contradiction completes the
%proof.
% $\hfill\blacksquare$

Combining Theorems \ref{main theorem3}, \ref{main theorem1} and
\ref{main theorem2} yields a positive answer to the question
\ref{question}. In the sequel, we determine
$R(\mathcal{P}^3_m,\mathcal{C}^3_n)$ when $n>m\geq 3$.

%%%%%%%%%%%%%%%%%%%%%%%%%%%%
\bigskip
\begin{lemma}\label{main theorem}
$R(\mathcal{P}^3_3,\mathcal{C}^3_4)=9.$
\end{lemma}

%%%%%%%%%%%%%%%%%%%%%%%%%%%%
%%%%%%%%%%%%%%%%%%%%%%%%%%%%
\noindent\textbf{Proof.} Using Theorem \ref{R(Pk3,Pk3)} we have
$R(\mathcal{C}^3_4,\mathcal{C}^3_4)=9$. On the other hand
$R(\mathcal{P}^3_3,\mathcal{C}^3_4)\leq
R(\mathcal{C}^3_4,\mathcal{C}^3_4)$.
 $\hfill\blacksquare$
%%%%%%%%%%%%%%%%%%%%%%%%%%%%%%%%%%%%%%%%%%%%%%%%%%%%%%%%%
%%%%%%%%%%%%%%%%%%%%%%%%%%%%%%%%%%%%%%%%%%%%%%%%%%%%%%%%%%%%%%%%%%%%%%%%%%%%%%%%%%%%%%%%%%%%%%%%%%%%%%

 \bigskip
\begin{theorem}\label{main theorem4}
%We have
%\begin{eqnarray*}
%R(\mathcal{P}^3_n,\mathcal{C}^3_m)= \left\lbrace
%\begin{array}{ll}
%2n+\Big\lfloor\frac{m+1}{2}\Big\rfloor &\mbox{if~}  n\geq m\geq 3,\vspace{.5 cm}\\
%2m+\Big\lfloor\frac{n-1}{2}\Big\rfloor &\mbox{if~} m> n\geq 3.
%\end{array}
%\right.\vspace{.2 cm}
%\end{eqnarray*}
For every $n>m\geq 3$,
$$R(\mathcal{P}^3_m,\mathcal{C}^3_n)=2n+\Big\lfloor\frac{m-1}{2}\Big\rfloor.$$
\end{theorem}

%%%%%%%%%%%%%%%%%%%%%%%%%%%%
\noindent\textbf{Proof. } We prove the theorem by induction on
$m+n$.  By Lemma \ref{main theorem} the base case is trivial.
Suppose to the contrary that
$\mathcal{H}=\mathcal{K}^3_{2n+\lfloor\frac{m-1}{2}\rfloor}$ is
2-edge colored red and blue with  no a red copy of
$\mathcal{P}^3_m$  and no a blue copy of $\mathcal{C}^3_n$ in
$\mathcal{H}$.
 Consider the following cases.

\bigskip
\noindent \textbf{Case 1. } $n=m+1$

\medskip
By Theorem \ref{main theorem2} we have
$R(\mathcal{P}^3_{m},\mathcal{C}^3_{m})=
2m+\Big\lfloor\frac{m+1}{2}\Big\rfloor< 2(m+1)+
\Big\lfloor\frac{m-1}{2}\Big\rfloor$. Since there is no red copy
of $\mathcal{P}^3_{m}$,  we have a blue $\mathcal{C}^3_{m}$. Now,
by using Lemma \ref{cn-1 implies cm} there is a red copy of
$\mathcal{P}^3_{m}$,
 a contradiction.

\bigskip
\noindent \textbf{Case 2. }$n> m+1$.

\medskip

By the induction hypothesis
 $R(\mathcal{P}^3_{m},\mathcal{C}^3_{n-1})=
2(n-1)+\Big\lfloor\frac{m-1}{2}\Big\rfloor<
2n+\Big\lfloor\frac{m-1}{2}\Big\rfloor$. Since there is no red
copy of $\mathcal{P}^3_m$, we have a copy of $\mathcal{C}^3_{n-1}$
in $\mathcal{H}_{blue}$.
 So using  Lemma \ref{cn-1 implies cm} we have a red copy
 of  $\mathcal{P}^3_{m}$, a contradiction.
 $\hfill\blacksquare$

%%%%%%%%%%%%%%%%%%%%%%%%%%%%%%%%%%%%%%%%%%%%%%%%%%%%%%%%%%%%%%%%%%%%%%%%%%%%%%%%%%%%%%%%%%%%%%%%%%%%%%%%%%

%\vspace{0.5 cm} It is worth noting that  we can also conclude
%$R(\mathcal{C}^3_n,\mathcal{P}^3_m)\leq
%n+\lfloor\frac{m-1}{2}\rfloor$, for $n
%> m\geq 3$. To see that assume
%$\mathcal{K}^3_{2n+\lfloor\frac{m-1}{2}\rfloor}$ be 2-edge
%colored red and blue with  no a red copy of $\mathcal{C}^3_n$
%and a  blue copy of $\mathcal{P}^3_m$. Since
%$R(\mathcal{P}^3_{n-1},\mathcal{P}^3_m)=
%2(n-1)+\lfloor\frac{m+1}{2}\rfloor <
%2n+\lfloor\frac{m-1}{2}\rfloor$, there is a copy of
%$\mathcal{P}^3_{n-1}$ in $\mathcal{H}_{red}$. So using Lemma
%\ref{pn-1 implies cm and pm} we have a copy of $\mathcal{P}^3_m$
%in $\mathcal{H}_{blue}$, a contradiction. So we have the
%following.

%%%%%%%%%%%%%%%%%%%%%%%%%%%%

% \bigskip
%\begin{theorem}\label{main theorem4}
%We have
%\begin{eqnarray*}
%R(\mathcal{P}^3_n,\mathcal{C}^3_m)= \left\lbrace
%\begin{array}{ll}
%2n+\Big\lfloor\frac{m+1}{2}\Big\rfloor &\mbox{if~}  n\geq m\geq 3,\vspace{.5 cm}\\
%2m+\Big\lfloor\frac{n-1}{2}\Big\rfloor &\mbox{if~} m> n\geq 3.
%\end{array}
%\right.\vspace{.2 cm}
%\end{eqnarray*}
%For every $m>n\geq 3$,
%$$R(\mathcal{P}^3_n,\mathcal{C}^3_m)=2m+\Big\lfloor\frac{n-1}{2}\Big\rfloor.$$
%\end{theorem}

%%%%%%%%%%%%%%%%%%%%%%%%%%%%

%%%%%%%%%%%%%%%%%%%%%%%%%%%%%%%%%%%%%%%%%%%%%%%%%%%%%%%%%%%%%%%%%%%%%%%%%%%%%%%%%%%%%%%%%%%%%%%%%%%%%%%%%%
%%%%%%%%%%%%%%%%%%%%%%%%%%%%%%%%%%%%%%%%%%%%%%%%%%%%%%%%%%%%%%%%%%%%%%%%%%%%%%%%%%%%%%%%%%%%%%%%%%%%%%%%%%%
%%%%%%%%%%%%%%%%%%%%%%%%%%%%%%%%%%%%%%%%%%%%%%%%%%%%%%%%%%%%%%%%%%%%%%%%%%%%%%%%%%%%%%%%%%%%%%%%%%%%%%%%%%%

\end{document}